\documentstyle[12pt,amscd,amssymb]{amsart}

\def\sw#1{{\sb{(#1)}}}
\def\tens{\mathop{\otimes}}

\def\<{{\langle}}
\def\>{{\rangle}}

\def\eps{\epsilon}

\def\note#1{{}}

\def\note#1{}
\def\M{{\bf M}}

\def\cC{{\cal C}}
\def\cA{{\cal A}}

\def\Label{\label}
\begin{document}

\baselineskip 22pt

\newtheorem{proposition}{Proposition}[section]
\newtheorem{lemma}[proposition]{Lemma}
\newtheorem{corollary}[proposition]{Corollary}
\newtheorem{theorem}[proposition]{Theorem}

\theoremstyle{definition}
\newtheorem{definition}[proposition]{Definition}
\newtheorem{example}[proposition]{Example}

\theoremstyle{remark}
\newtheorem{remark}[proposition]{Remark}

\newcommand{\Section}{\setcounter{definition}{0}\section}
\renewcommand{\theequation}{\thesection.\arabic{equation}}
\newcounter{c}
\renewcommand{\[}{\setcounter{c}{1}$$}
\newcommand{\etyk}[1]{\vspace{-7.4mm}$$\begin{equation}\Label{#1}
\addtocounter{c}{1}}
\renewcommand{\]}{\ifnum \value{c}=1 $$\else \end{equation}\fi}

\title[Induction functors and corings]{The structure of
corings.\\~\\ {\em Induction functors, Maschke-type theorem, and Frobenius and
Galois-type properties.}}
\author{Tomasz Brzezi\'nski}
\address{Department of Mathematics, University of Wales Swansea,
Singleton Park, Swansea SA2 8PP, U.K.}
\address{Department of Theoretical Physics, 
University of
\L\'od\'z, Pomorska 149/153, 90-236, \L\'od\'z, Poland}
\email{T.Brzezinski@@swansea.ac.uk}
\thanks{EPSRC Advanced Research Fellow.}
\begin{abstract}
Given a ring $A$ and an $A$-coring $\cC$ we study when the forgetful
functor from the category of right $\cC$-comodules to the category of
right $A$-modules and its right adjoint
$-\otimes_A\cC$ are separable. We then proceed to study when the
induction functor $-\otimes_A\cC$  is also the left adjoint of the
forgetful functor. This question is closely related to the problem when
$A\to {}_A{\rm Hom}(\cC,A)$ is a Frobenius extension.
 We introduce the notion of a Galois
coring and analyse when the tensor functor over the subring of
$A$ fixed under the coaction of $\cC$ is an equivalence. We also
comment on possible dualisation of the notion of a coring.
\end{abstract}
\maketitle
\section{Introduction}
A coring is a generalisation of a coalgebra introduced by M.\ Sweedler
in \cite{Swe:pre}. It has been recently pointed out by M.\ Takeuchi
\cite{Tak:com} that new examples of corings can be provided by 
{\em entwining structures} introduced in \cite{BrzMa:coa} in the context
of gauge theory on non-commutative spaces. Entwining
structures and modules associated to them generalise the
notion of Doi-Koppinen Hopf modules introduced in \cite{Doi:uni}
\cite{Kop:var}. Various structure theorems concerning Doi-Koppinen
modules can be formulated more generally in terms of entwined modules
(see recent paper \cite{CaeDeG:mod} or \cite{Brz:coa} 
for an exhaustive list of
references). In the present paper we argue that many of those structure
theorems are in fact special cases of structure theorems for the
category of comodules of a coring.

We begin in Section~2 by recalling the definition of a coring and
listing examples of corings coming from entwining structures and their
recent generalisation \cite{CaeDeG:mod} motivated by Doi-Koppinen
modules for a weak Hopf algebra introduced in \cite{Boh:Doi}. In
Section~3 we study when the forgetful functor from the category of right
comodules of an $A$-coring $\cC$ to the category of right $A$-modules is
separable. This turns out to be equivalent to coseparability of $\cC$
and implies a Maschke-type theorem for a coring. We also study
when the induction functor $-\otimes_A\cC$ which is the right adjoint of the
forgetful functor above is separable. In Section~4 we derive conditions
for the induction functor to be the left adjoint of the forgetful functor
as well. Put differently, since a functor which has the same left and right adjoint is
known as a {\em Frobenius functor} \cite{CaeMil:Doi}, we derive in Section~4 conditions
 for the forgetful functor (or, equivalently, the induction functor)
to be a Frobenius functor. 
This problem  is
closely  related to the question when $A\to {}_A{\rm Hom}(\cC, A)$ is a Frobenius
extension. Such a question is particularly relevant in the present
context as it is known that there is a close relationship between
Frobenius extensions $A\to R$ and certain $A$-coring structures on $R$
\cite{CaeIon:fro} \cite{Kad:new}. Next in Section~5 we analyse the
situation in which the ring $A$ is itself a right $\cC$-comodule. 
We define
the subring of coinvariants $B$ of $A$ and study the induction 
functor $-\otimes_B
A$ from the category of right $B$-modules to the category of right
$\cC$-comodules. Motivated by coalgebra-Galois extensions we introduce 
the notion of a Galois coring and show that if $\cC$ is a flat left
$A$-module then $-\otimes_B A$ is an
equivalence if and only if $A$ is a faithfully flat left $B$-module and
$\cC$ is Galois. Finally, in Section~6 we define  a
$C$-algebroid or a $C$-ring 
as a dualisation of $A$-coring and gather some comments on
dualisation of the results of the previous sections to the case of 
$C$-algebroids.  

We use the following conventions. For an object $V$ in a category, 
the identity morphism $V\to V$ is
 denoted by $V$. All rings in this paper have 1, a ring map is 
assumed to respect 1, and unless stated otherwise all modules over a ring 
are assumed to be unital.
For a ring $A$,  $\M_A$ (resp.\ ${}_A\M$) denotes 
the category of right (resp.\ left) $A$-modules. The morphisms in this
category are denoted by ${\rm Hom}_A(M,N)$ (resp.\ ${}_A{\rm Hom}(M,N)$).
 The action of $A$ is denoted by a
dot between elements. 
If $A,B$ are rings and $M,N$ are $(A,B)$-bimodules then ${}_A{\rm
Hom}_B(M,N)$ denotes the set of $(A,B)$-bimodule maps. For an
$(A,A)$-bimodule $M$, $M^A$ denotes the subbimodule of invariants, i.e.,
$M^A = \{m\in M\;|\; \forall a\in A,\; a\cdot m = m\cdot a\}$.

Throughout the paper $k$ denotes a commutative ring with unit. 
We  assume that all
the algebras are over $k$ and  unital, and  coalgebras are over $k$
and  counital. Unadorned tensor product is over $k$.
For a $k$-algebra $A$ we use $\mu$ to denote the product as a map and 
$1$ to denote  unit both as an element of $A$ and as a map $k\to
A$, $\alpha\to \alpha 1$. 

For a
$k$-coalgebra $C$ we use $\Delta$ to denote the coproduct and $\eps$ to
denote the counit. Notation for comodules is similar to that for modules
but with subscripts replaced by superscripts, i.e. $\M^C$ is the
category of right $C$-comodules, $\rho^M$ is a right coaction, ${}^C\M$
is the category of left $C$-comodules and ${}^M\rho$ is a left coaction
 etc. We
use the Sweedler notation for coproducts and coactions, i.e. $\Delta(c)
= c\sw 1\otimes c\sw 2$, $\rho^M(m) = m\sw 0\otimes m\sw 1$ (summation
understood). We use similar notation for corings.

\section{New examples of corings}
Let $A$ be a ring. Recall from \cite{Swe:pre} that an 
{\em $A$-coring} is an $(A,A)$-bimodule
$\cC$ together with $(A,A)$-bimodule maps $\Delta_\cC:\cC\to
\cC\otimes_A\cC$ called a coproduct and $\eps_\cC:\cC\to A$ called a 
counit, such that
$$
(\Delta_\cC\otimes_A\cC)\circ\Delta_\cC = (\cC\otimes_A\Delta_\cC)\circ
\Delta_\cC, \quad (\eps_\cC\otimes_A\cC)\circ\Delta_\cC =
(\cC\otimes_A\eps_\cC)\circ \Delta_\cC = \cC.
$$
 Given an $A$-coring $\cC$ 
a {\em right $\cC$-comodule} is a right $A$-module $M$ together with a right
$A$-module map $\rho^M:M\to M\otimes_A\cC$, called a {\em coaction}
  such that
$$
(\rho^M\otimes_A\cC)\circ\rho^M = (M\otimes_A\Delta_\cC)\circ
\rho^M, \quad 
(M\otimes_A\eps_\cC)\circ \rho^M = M.
$$
A map between  right $\cC$-comodules is a right $A$-module map
respecting the coactions, i.e., $f:M\to N$ satisfies
$\rho^N\circ f = (f\otimes_A\cC)\circ\rho^M$. The category of
right $\cC$-comodules is denoted by $\M^\cC$, and their morphisms by
${\rm Hom}^\cC(-,-)$.
\begin{example}(\cite[1.2~Example]{Swe:pre})
Let $B\to A$ be a ring extension. Then $\cC = A\otimes_B A$ is an
$A$-coring with the coproduct $\Delta_\cC : \cC\to \cC\otimes_A\cC\cong
A\otimes_BA\otimes_BA$, $a\otimes_B a'\mapsto a\otimes_B 1\otimes_B a'$ and the
counit $\eps_\cC(a\otimes_B a') = aa'$. $\cC$ is called the {\em
canonical $A$-coring} associated to the extension $B\to A$.
\label{can.coring}
\end{example}

Interesting examples of corings come from entwining structures. 
Recall from \cite{BrzMa:coa} that
an {\em  entwining
structure} (over $k$) is a 
triple $(A,C)_\psi$ consisting of a $k$-algebra $A$, a $k$-coalgebra
$C$ and a $k$-module map $\psi: C\tens A\to A\tens C$ satisfying
\begin{equation}
\psi\circ(C\tens \mu) = (\mu\tens C)\circ (A\tens\psi)\circ(\psi\tens A),
\label{diag.A}
\end{equation}
\begin{equation}
(A\tens\Delta)\circ\psi = (\psi\tens
C)\circ(C\tens\psi)\circ(\Delta\tens A),
\label{diag.B}
\end{equation}
\begin{equation}
\psi\circ (C\tens 1) = 1\tens C, \qquad (A\tens \eps)\circ\psi =
\eps\tens A.
\label{unit}
\end{equation}
For $(A,C)_\psi$ we  use the
notation $\psi(c\otimes a) = a_\alpha \otimes c^\alpha$ (summation
over a Greek index understood), for all $a\in
A$, $c\in C$. Various examples of entwining structures can be found in
\cite{Brz:mod}. Given an entwining structure $(A,C)_\psi$, an (entwined) 
{\em $(A,C)_\psi$-module}
is a right $A$-module, right $C$-comodule $M$ 
such that 
\begin{equation}
\rho^M\circ\rho_M = (\rho_M\tens C)\circ(M\tens\psi)\circ(\rho ^M\otimes
A), 
\label{entw.mod}
\end{equation}
where $\rho_M$ is the action of $A$ on $M$. Explicitly: 
$
\rho^M(m\cdot a) = m\sw 0\cdot a_\alpha\tens m\sw 1^\alpha$, $ \forall
a\in A, m\in M$. $A\otimes C$ is an example of entwined module with
the coaction $A\otimes \Delta$ and the action $(a'\otimes c)\cdot a =
a'\psi(c\otimes a)$, for all $a,a'\in A, c\in C$. The unifying
Doi-Koppinen modules \cite{Doi:uni} \cite{Kop:var} are another example of
entwined modules  (cf.\ \cite[Example~3.1(3)]{Brz:mod}).

A morphism of $(A,C)_\psi$-modules is a right $A$-module map
which is also a right $C$-comodule map. The category of
$(A,C)_\psi$-modules is denoted by $\M_A^C(\psi)$.  

The present paper is motivated by the following
observation ascribed to M.\ Takeuchi \cite{Tak:com}
\begin{proposition}
For an entwining structure $(A,C)_\psi$, view $A\otimes C$ as an
$(A,A)$-bimodule with the left action $a\cdot(a'\otimes c)= aa'\otimes c$
and the right action $(a'\otimes c)\cdot a =
a'\psi(c\otimes a)$, for all $a,a'\in A, c\in C$.  Then 
$\cC = A\otimes C$ is an
$A$-coring with the coproduct
 $\Delta_\cC : \cC\to \cC\otimes_A\cC\cong
A\otimes C\otimes C$, $\Delta_\cC = A\otimes \Delta$,  and the counit 
$\eps_\cC = A\otimes\eps$.

Conversely if $A\otimes C$ is an $A$-coring with the coproduct, 
counit and left $A$-module structure above, 
then $(A,C)_\psi$ is an entwining structure, where $\psi:
c\otimes a \mapsto (1\otimes c)\cdot a$.

Under this bijective correspondence $\M^\cC = \M_A^C(\psi)$.
\label{prop.tak}
\end{proposition}
\begin{proof}
One can easily check that given an entwining structure $(A,C)_\psi$ 
one has the coring $A\otimes C$ as described in the proposition.
Conversely, let $A\otimes C$ be an $A$-coring with structure maps given
in the proposition.  The
properties of the right $A$-action imply equation (\ref{diag.A}) and 
the first of equations (\ref{unit}) 
required for the entwining map $\psi$. The remaining two 
conditions follow 
from the
facts that $A\otimes\Delta$ and $A\otimes\eps$ are right $A$-module maps
respectively.
\end{proof}

Recall from \cite{BrzHaj:coa} that given a coalgebra $C$, an extension of algebras $B\subset A$
is $C$-Galois if $A$ is a right $C$-comodule with coaction $\rho^A$, $B
= A^{coC}=\{ b\in A\; |\; \forall a\in A, \; \rho^A(ba) = b\rho^A(a)\}$,
and the canonical left $A$-module, right
$C$-comodule map $can: A\otimes_B A\to A\otimes C$, $a\otimes_B
a'\mapsto a\rho^A(a')$ is bijective. In this case, since $A\otimes_BA$
is an $A$-coring by Example~\ref{can.coring}, the canonical map $can$ induces
an $A$-coring structure on $A\otimes C$. The corresponding entwining map
$\psi$  computed from Proposition~\ref{prop.tak} comes out as
$\psi(c\otimes a) = can(can^{-1}(1\otimes c)a)$, for all $a\in A$ and
$c\in C$ and thus coincides with the canonical
entwining structure associated to a $C$-Galois extension $B\subset A$
constructed in \cite[Theorem~2.7]{BrzHaj:coa}.

In the context of $C$-Galois extensions one can also mention the
following example of a coring inspired by \cite[3.5~Theorem]{Sch:big}. 
Given 
a $C$-Galois extension $B\subset A$ define 
a $(B,B)$-bimodule 
$$
\cC = \{ \sum_i a^i\otimes \overline{a}^i\in
A\otimes A\; |\; \sum_i a^i\sw 0\otimes can^{-1}(1\otimes a^i\sw 1)
\overline{a}^i = \sum_i a^i\otimes \overline{a}^i\otimes_B 1\}.
$$
If $A$ is 
faithfully flat as both $k$ and $B$-bimodule then $\cC$ is a
$B$-coring with the coproduct and counit
$$
\Delta_\cC (\sum_i a^i\otimes \overline{a}^i) = \sum_i a^i\sw 0\otimes
can^{-1}(1\otimes a^i\sw 1) \otimes \overline{a}^i, \qquad \eps_\cC (\sum_i a^i\otimes
\overline{a}^i) = \sum_i a^i\overline{a}^i.
$$
The fact that $\cC$ is a coring can be verified by direct computation
which uses properties of the canonical map $can$.

A generalisation of the notion of an entwining structure is possible by
replacing equations (\ref{unit}) by  weaker
conditions, for all $a\in A$, $c\in C$,
\begin{equation}
 a_\alpha\eps(c^\alpha) = 1_\alpha a\eps(c^\alpha), \quad 1_\alpha
\eps(c\sw 1^\alpha) \otimes c\sw 2 = 1_\alpha\otimes c^\alpha, 
\label{w.unit}
\end{equation}
where $\psi(c\otimes a) = a_\alpha\otimes c^\alpha$. 
Such a weakened entwining structure is termed a {\em weak entwining
structure} \cite{CaeDeG:mod} as it includes examples coming from weak 
Hopf algebras \cite{Boh:Doi}\cite{Boh:wea}, and is denoted by
$(A,C,\psi)$. 
Given a weak entwining structure $(A,C,\psi)$ one can
define the category of weak entwined modules as right $A$-modules, right
$C$-comodules such that eq.~(\ref{entw.mod}) holds. We denote this
category by $\M_A^C(\psi)$ too.
\begin{proposition}
Let $(A,C,\psi)$ be a weak entwining structure. Let $p :A\otimes C\to
A\otimes C$, $p = (\mu\otimes C)\circ(A\otimes \psi)\circ(A\otimes
C\otimes 1)$, and $\cC = {\rm Im}p = 
\{\sum_ia_i1_\alpha\otimes c_i^\alpha \; |\; \sum_ia_i\otimes c_i\in
A\otimes C\}$. Then $p\circ p =p$ and 

(1) $\cC$ is an
$(A,A)$-bimodule with the left action $a'\cdot(a1_\alpha\otimes
c^\alpha) = a'a1_\alpha\otimes c^\alpha$ and the right action
$(a'1_\alpha \otimes c^\alpha)\cdot a = a'1_\alpha a_\beta \otimes 
c^{\alpha\beta}= a'a_\alpha\otimes c^\alpha$ (this last equality follows
from equation (\ref{diag.A})). 

(2) $\cC$ is an $A$-coring with coproduct $\Delta_\cC = (A\otimes
\Delta)\mid_\cC$ and counit $\eps_\cC = (A\otimes \eps)\mid_\cC$.

(3) $\M_A^C(\psi) = \M^\cC$
\label{weak}
\end{proposition} 
\begin{proof} For any $a\otimes c\in A\otimes C$, using the definition of
$p$ and equation (\ref{diag.A}) we have
$$
p\circ p(a\otimes c) = p(a1_\alpha\otimes c^\alpha) =a1_\alpha
1_\beta \otimes c^{\alpha\beta} = a1_\alpha\otimes c^\alpha = 
p(a\otimes c).
$$
Therefore $p$ is a projection as claimed.

(1) Clearly $\cC$ is a left $A$-module. The fact that it is a right
$A$-module can be checked directly using equations (\ref{diag.A}) and
(\ref{w.unit}). 

(2) First we need to show that ${\rm Im}\Delta_\cC \subset
\cC\otimes_A\cC$. Since $\Delta_\cC$ is $k$-linear and, evidently,
 left $A$-linear, it suffices to take an element of $\cC$ of the form
$e= 1_\alpha\otimes c^\alpha$ and compute
\begin{eqnarray*}
\Delta_\cC(e) &=& 1_\alpha\otimes c^\alpha\sw
1\otimes c^\alpha\sw 2 = 1_{\alpha\beta}\otimes c\sw 1^\beta\otimes
c\sw 2^\alpha = 1_\gamma 1_{\alpha\beta} \otimes c\sw 1^{\gamma\beta}
\otimes c\sw 2^\alpha\\
&=& (1_\gamma\otimes c\sw 1^\gamma)\cdot 1_\alpha\otimes c\sw 2^\alpha =
(1_\gamma \otimes c\sw 1^\gamma)\otimes_A(1_\alpha\otimes c\sw
1^\alpha)\in \cC\otimes_A\cC.
\end{eqnarray*}
We used equations (\ref{diag.B}), (\ref{diag.A}) and the fact that $1$
is a unit in $A$ to derive second and third equalities, and then
definition of the right action of $A$ on $\cC$ to obtain the fourth one.

To prove that $\Delta_\cC$ is a right $A$-module map it is enough to
consider any $e\in\cC$ of the above form, use the same equations as for
the preceding calculation, and compute for all $a\in A$,
\begin{eqnarray*}
\Delta_\cC(e\cdot a) &=& a_\alpha\otimes c^\alpha\sw 1\otimes
c^\alpha\sw 2 = a_{\alpha\beta}\otimes c\sw 1^\beta\otimes
c\sw 2^\alpha  = (1_\gamma\otimes c\sw
1^\gamma)\otimes_A(a_\alpha\otimes c\sw 2^\alpha)\\
&=& (1_\gamma\otimes c\sw1^\gamma)\otimes_A(1_\alpha\otimes c\sw
2^\alpha)\cdot a =\Delta_\cC(e)\cdot a.
\end{eqnarray*}
The coassociativity of $\Delta_\cC$ follows from the coassociativity of
the coproduct $\Delta$. Now, it is clear from the definition that
$\eps_\cC$ is a left $A$-module map. To prove that it is a right
$A$-module map as well take $e\in\cC$ as above and compute for all $a\in
A$, 
$$
\eps_\cC(e\cdot a) =a_\alpha\eps(c^\alpha) = 1_\alpha a\eps(c^\alpha)  =
\eps_\cC(e)a,
$$
where the first of equations (\ref{w.unit}) was used to obtain the
second equality. Again, directly from the construction of $\eps_\cC$ it
follows that $(\cC\otimes_A\eps_\cC)\circ\Delta_\cC = \cC$. On the other
hand using the definitions of $\Delta_\cC$ and $\eps_\cC$, the second of equations
(\ref{w.unit}) and equation (\ref{diag.A}) we obtain for any $e\in\cC$ of the above form
$$
(\eps_\cC\otimes_A\cC)\circ\Delta_\cC(e) = \eps_\cC(1_\alpha\otimes c\sw
1^\alpha)\cdot (1_\beta\otimes c\sw 2^\beta) = 1_\alpha1_\beta\eps(c\sw
1^\alpha)\otimes c\sw 2^\beta = 1_\alpha1_\beta\otimes c^{\alpha\beta} =
e.
$$
This completes the proof that $\cC$ is an $A$-coring.

(3) If $M\in\M_A^C(\psi)$ with the right $C$-coaction $m\mapsto m\sw
0\otimes m\sw 1$, then $M\in\M^\cC$ with the coaction $\rho^M(m) = m\sw
0\otimes _A(1_\alpha\otimes m\sw 1^\alpha) = m\sw 0\cdot 1_\alpha
\otimes m\sw
1^\alpha = m\sw 0\otimes m\sw 1$ since $M\in\M_A^C(\psi)$. 
Conversely, if $M\in\M^\cC$
then the coaction $\rho^M :M\to M\otimes_A\cC \subset M\otimes C$
provides $M$ also with the weak entwined module structure.
\end{proof}

An example of a weak entwining structure can be obtained by the
following construction.
\begin{example}
Let $C$ be a coalgebra, $A$ an algebra and a right $C$-comodule with the
coaction $\rho^A$. Let
$B=A^{coC}=\{ b\in A\; |\; \forall a\in A, \; \rho^A(ba) = b\rho^A(a)\}$,
and let  $can:A\otimes_BA\to A\otimes C$, $a\otimes_B a'\mapsto
a\rho^A(a')$. View $A\otimes_B A$ as a left $A$-module via
$\mu\otimes_B A$ and a right $C$-comodule via $A\otimes_B\rho^A$. View
$A\otimes C$ as a left $A$-module via $\mu\otimes C$ and as a right
$C$-comodule via $A\otimes\Delta$. Now suppose that  $can$ is a split
monomorphism in the category of left $A$-modules and right
$C$-comodules, i.e., there exists left
$A$-module, right $C$-comodule map $\sigma:A\otimes C\to A\otimes_B A$
such that $\sigma\circ can =A\otimes_B A$. Let 
 $\tau :C\to A\otimes_BA$, $c\mapsto \sigma(1\otimes c)$. Define
$$
\psi:C\otimes A\to A\otimes C, \quad \psi
=can\circ(A\otimes_B\mu)\circ(\tau\otimes A).
$$
Then $(A,C,\psi)$ is a weak entwining structure. The extension of
algebras $B\subset A$ is called a {\em weak $C$-Galois extension}.
\label{ex.weak}
\end{example}
Example~\ref{ex.weak} can be proven by a direct verification of the
axioms in  a way similar to
the proof of \cite[Theorem~2.7]{BrzHaj:coa}. Similarly one can prove 
that in this case $A$ is a weak entwined $(A,C,\psi)$-module via 
$\rho^A$ and multiplication in
$A$, and that $(A,C,\psi)$ is unique with this respect.

The fact that given a (weak) entwining structure $(A,C)_\psi$, 
there is an associated
$A$-coring $\cC$ explains various properties of entwining structures. 
For example, 
the existence of the $\psi$-twisted convolution product on ${\rm
Hom}(C,A)$ given by $f*_\psi g(c) = f(c\sw 2)_\alpha g(c\sw 1^\alpha)$ is
the consequence of the product on ${}_A{\rm Hom}(\cC, A)$
defined in \cite[3.2~Proposition~(a)]{Swe:pre}. We deal with this product
in Section~4. In the remaining part of
this paper we show that some properties of entwined modules can be
derived from properties of comodules of a coring. 

\section{Separable functors for a coring and a Maschke-type theorem}
In this section, for an $A$-coring $\cC$ we study when two functors, the
forgetful functor $\M^\cC\to\M_A$ and the
induction  functor 
$-\otimes_A\cC :\M_A\to \M^\cC$, $M\mapsto M\otimes_A\cC$ are separable. The right 
$\cC$-comodule structure of $M\otimes_A\cC$ is given by
$M\otimes_A\Delta_\cC$. Recall from \cite{Nas:sep} that
a covariant functor $F :\ {\bf C}\to {\bf D}$ is {\em separable} if the
natural transformation ${\rm Hom}_{\bf C}(-,-)\to {\rm Hom}_{\bf
D}(F(-), F(-))$ splits. We start by observing the following fact 
\cite[Proposition~3.1]{Guz:coi}.
\begin{lemma}
The functor $-\otimes_A\cC$ is the right adjoint of the forgetful
functor $F:\M^\cC\to\M_A$.
\label{lemma.adjoint}
\end{lemma}
\begin{proof}
For any $M\in \M_A$ define, $\Psi_M :M\otimes_A\cC \to M$, $\Psi_M =
M\otimes_A\eps_\cC$, while for any $N\in \M^\cC$ define $\Phi_N: N\to
N\otimes_A\cC$, $\Phi_N = \rho^N$. Clearly both $\Psi_M$ and $\Phi_N$
are natural in $M$ and $N$ respectively. Furthermore, because $\rho^N$
is a coaction $\Psi_N\circ\Phi_N = N$, for all $N\in\M^\cC$. It remains
to be shown that for all right $A$-modules $M$,
$(\Psi_M\otimes_A\cC)\circ\Phi_{M\otimes_A\cC} = M\otimes_A\cC$. Take
any $m\in M$ and $c\in \cC$ and compute
$$
(\Psi_M\otimes_A\cC)\circ\Phi_{M\otimes_A\cC}(m\otimes_Ac) =
m\otimes_A\eps_\cC(c\sw 1)\cdot c\sw 2 = m\otimes_Ac.
$$
Thus $\Psi$
and $\Phi$ are the required adjunctions.
\end{proof}

In particular if $\cC$ is taken as in Proposition~\ref{weak},
Lemma~\ref{lemma.adjoint} implies \cite[Proposition~2.2]{CaeDeG:mod}. We
now study when the functors in Lemma~\ref{lemma.adjoint} are separable.
 Since we are dealing with a pair of
adjoint functors, the following characterisation of separable
functors, obtained in \cite{Raf:sep} \cite{Rio:cat}, is of  great
importance

\begin{theorem}
Let $G:\ {\bf D}\to {\bf C}$ be the right adjoint of $F:\ {\bf C}\to {\bf
D}$ with adjunctions $\Phi:\ {\bf C}\to GF$ and
$\Psi:\ FG\to {\bf D}$. Then

(1) $F$ is separable if and only if $\Phi$ splits, i.e., for all objects
$C\in{\bf C}$ there exists a morphism $\nu_C \in{\rm Mor}_{\bf C}
(GF(C), C)$ such that $\nu_C\circ \Phi_C = C$ and for all $f\in {\rm
Mor}_{\bf C} (C, \tilde{C})$, $\nu_{\tilde{C}}\circ GF(f) = f\circ\nu_C$.

(2) $G$ is separable if and only if $\Psi$ cosplits, i.e., for all
objects $D\in {\bf D}$ there exists a morphism $\nu_D \in{\rm Mor}_{\bf D}
(D, FG(D))$ such that $\Psi_D\circ \nu_D = D$ and for all $f\in {\rm
Mor}_{\bf D} (D, \tilde{D})$, $\nu_{\tilde{D}}\circ f = FG(f)\circ\nu_D$.
\label{thm.sep}
\end{theorem}

First we analyse when $-\otimes_A\cC$ is a separable functor.
\begin{theorem}
Let $\cC$ be an $A$-coring. Then the functor $-\otimes_A\cC$ is
separable if and only if there exists an invariant $e\in \cC^A$ such
that $\eps_\cC (e) = 1$.
\label{thm.sep1}
\end{theorem}
\begin{proof} ``$\Rightarrow$". Suppose that $-\otimes_A\cC$ is
separable, and let $\nu$ be split by $\Psi$, the latter defined in
Lemma~\ref{lemma.adjoint}. Since $A$ is a
right $A$-module we can take $\nu_A\in {\rm Hom}_A(A, A\otimes_A\cC) \cong {\rm Hom}_A(A,
\cC)$ and define $e=\nu_A(1)\in \cC$. Since $\nu_A$ is split
 by $\Psi_A=A\otimes_A\eps_\cC = \eps_\cC$, 
we have $1 = \Psi_A\circ\nu_A(1) = \eps_\cC\circ\nu_A(1) =
\eps_\cC(e)$. Now for any $a\in A$ define the morphism $f_a\in {\rm Hom}_A(A,
A)$ by $f_a(a') = aa'$. Since $\nu_A$ is natural we have:
$$
\nu_A(aa') =\nu_A\circ f_a(a') = (f_a\otimes_A\cC)\circ\nu_A(a') =
a\cdot\nu_A(a').
$$
This implies, in particular, $\nu_A(a) = a\cdot\nu_A(1) = a\cdot e$. On
the other hand, $\nu_A$ is a right $A$-module morphism, thus we have,
$a\cdot e = \nu_A(a) = \nu_A(1)\cdot a = e\cdot a$, so that $e$ is
$A$-central as required.

``$\Leftarrow$". Suppose there exists $e\in\cC^A$ such that 
$\eps_\cC(e)=1$. For any $M\in \M_A$ define
$\nu_M :M\to M\otimes_A\cC$ by $\nu_M(m) = m\otimes_A e$. Since $e$ is
$A$-central we have for all $a\in A$, $m\in M$:
$$
\nu_M(m\cdot a) = m\cdot a\otimes_A e = m\otimes_A a\cdot e = m\otimes_A
e\cdot a = \nu_M(m)\cdot a,
$$
hence $\nu_M$ is a right $A$-module map. Furthermore
$\Psi_M\circ\nu_M(m) = m\cdot\eps_\cC(e) = m$, i.e., $\nu_M$ is split
by 
$\Psi_M$.  Finally for any $f\in {\rm Hom}_A(M,N)$ we have $\nu_N\circ
f(m) = f(m)\otimes_A e = (f\otimes_A\cC)\circ \nu_M(m)$, which means
that $\nu_M$ is natural in $M$.
\end{proof}

One of the motivations for the introduction of separable functors was
that in the case of a ring extension $B\to A$, the restriction of scalars
functor is separable if and only if the extension $B\to A$ is separable
\cite{Nas:sep}. Theorem~\ref{thm.sep} provides one with  
different characterisation of separable
extensions. 
\begin{corollary}
Let $B\to A$ be a ring extension and let $\cC$ 
 be the canonical
$A$-coring in Example~\ref{can.coring}. The extension $B\to A$ 
is separable if and only if the
functor $-\otimes_A\cC$ is separable.
\label{cor1}
\end{corollary}
\begin{proof}
Recall the a ring extension  $B\to A$ is separable iff there
exists an invariant  $e = \sum_ia_i\otimes_B a'_i\in (A\otimes_BA)^A$  
such that $\sum_ia_ia'_i = 1$ (cf.\ \cite[Definition~2]{HirSug:sem}). Since in the canonical coring
$A\otimes_B A$ the counit is given by the multiplication projected to
$A\otimes_B A$ the separability of $B\to A$ is equivalent to the
separability of $-\otimes_A\cC$ by the preceding theorem.
\end{proof}

In the case of a coring $\cC = A\otimes C$ corresponding to an entwining
structure $(A,C)_\psi$ via Proposition~\ref{prop.tak}, a normalised
$A$-central element $e\in \cC$ is simply an integral in $(A,C)_\psi$ in
the sense of \cite[Definition~3.1]{Brz:coa}. Thus Theorem~\ref{thm.sep1}
implies \cite[Theorem~3.3]{Brz:coa} (which itself is a generalisation of
\cite[Theorem~2.14]{CaeMil:sep} formulated for Doi-Koppinen modules) while Corollary~\ref{cor1}
implies
\cite[Proposition~3.5]{Brz:coa}. 
\begin{theorem}
For an $A$-coring  $\cC$, the forgetful functor $F:\M^\cC\to\M_A$
is separable if and only if there exists $\gamma\in {}_A{\rm
Hom}_A(\cC\otimes_A\cC, A)$ such that $\gamma\circ\Delta_\cC = \eps_\cC$
and
\begin{equation}
 c\sw 1\cdot\gamma(c\sw 2\otimes_A c') =
\gamma(c\otimes_A c'\sw 1)\cdot c'\sw 2, \quad \forall c,c'\in \cC.
\label{eq}
\end{equation}
\label{thm.sep2}
\end{theorem}
\begin{proof}
``$\Rightarrow$". Suppose $F$ is separable. Let $\Phi$ be  
the adjunction defined in Lemma~\ref{lemma.adjoint}
and let $\nu$ be the splitting of  $\Phi$.  
Since $\cC$ is a right $\cC$-comodule via
$\Delta_\cC$ there is a corresponding $\nu_\cC$ and we can define $\gamma =
\eps_\cC\circ\nu_\cC :\cC\otimes_A\cC\to A$. The map $\gamma$ is a right
$A$-module morphism as a composition of two such morphisms. Next for all
$a\in A$ consider $f_a\in {\rm Hom}^\cC(\cC,\cC)$ given by $f_a(c) =
a\cdot c$. The naturality of $\nu$ implies for all $c,c'\in \cC$, 
$\nu_\cC(f_a(c)\otimes_A c') = f_a\circ\nu_\cC(c\otimes_A c')$, i.e.,
$\nu_\cC(a\cdot c\otimes_A c') = a\cdot\nu_\cC(c\otimes_A c')$. Therefore
$\nu_\cC$ is a left $A$-module map, and, consequently, $\gamma$ is a
left $A$-module morphism as a composition of two such morphisms. Since
$\nu_\cC$ splits $\Phi_\cC = \Delta_\cC$ we have
$$
\gamma\circ\Delta_\cC = \eps_\cC\circ\nu_\cC\circ\Delta_\cC = \eps_\cC.
$$
Now, for all $c\in \cC$ consider the morphism $\ell_c \in {\rm Hom}^\cC
(\cC,\cC\otimes_A\cC)$, $\ell_c(c') = c\otimes_A c'$, and also the
morphism $\Delta_\cC \in {\rm Hom}^\cC
(\cC,\cC\otimes_A\cC)$. By the naturality of $\nu_\cC$ we have
$$
\ell_c\circ\nu_\cC = \nu_{\cC\otimes_A\cC}\circ (\ell_c\otimes_A\cC),
\qquad \Delta_\cC\circ\nu_\cC = \nu_{\cC\otimes_A\cC}\circ (\Delta_\cC
\otimes_A\cC).\eqno{(*)}
$$
From the first of equations ($*$) we have for all $c,c',c''\in\cC$,
$c\otimes_A \nu_\cC(c'\otimes_A c'') = \nu_{\cC\otimes_A\cC}(c\otimes_A
c'\otimes_A c'')$. Combining this expression with the second of
equations ($*$) we
obtain
$$
\Delta_\cC (\nu_\cC(c\otimes_A c')) = \nu_{\cC\otimes_A\cC}(c\sw
1\otimes_A 
c\sw 2\otimes_A c') = c\sw 1\otimes_A \nu_\cC(c\sw 2\otimes_A c').
$$
Applying $\cC\otimes_A\eps_\cC$ to this equality one obtains
$$
\nu_\cC(c\otimes_A c') = c\sw 1\cdot\gamma (c\sw 2\otimes_A c').
$$ 
On the other hand, $\nu_\cC$ is a right $\cC$-comodule map, so that 
$\Delta_\cC(\nu_\cC (c\otimes_A c')) = \nu_\cC(c\otimes_A c'\sw
1)\otimes_A 
c'\sw 2$. Applying $\eps_\cC \otimes_A\cC$ to this formula 
we obtain
$$
\nu_\cC(c\otimes_A c') = \gamma (c\otimes_A c'\sw 1)\cdot c'\sw 2.
$$
Thus $\nu_\cC$ can be expressed in terms of $\gamma$ 
in two different ways. Comparison gives eq.~(\ref{eq}) as required.

``$\Leftarrow$". Suppose there exists $\gamma$ as in the theorem. Then
for all $M\in \M^\cC$ define an additive map
$\nu_M:M\otimes_A\cC\to M$, $m\otimes_A c\mapsto m\sw
0\cdot\gamma(m\sw 1\otimes_A c)$.
Clearly $\nu_M$ is a right $A$-module map since $\gamma$ is such a map.
Furthermore, for any $m\in M$, $c\in \cC$ we have
\begin{eqnarray*}
\nu_M(m\otimes_A c\sw 1)\otimes_A c\sw 2 &=& m\sw 0\cdot\gamma(m\sw 1
\otimes_A
c\sw 1)\otimes_A c\sw 2 \\
& = & m\sw 0\otimes_A \gamma(m\sw 1\otimes_A c\sw
1)\cdot c\sw 2 \\
&=&m\sw 0 \otimes_A m\sw 1\cdot\gamma(m\sw 2\otimes _Ac)\\
& = &
\nu_M(m\otimes_A 
c)\sw 0\otimes_A \nu_M(m\otimes_A c)\sw 1,
\end{eqnarray*}
where we used eq.~(\ref{eq}). 
This means that $\nu_M$ is a morphism in $\M^\cC$. Furthermore, take any
$f\in{\rm Hom}^\cC(M,N)$. Then for all $m\in M$, $c\in \cC$ we have
$$
\nu_M(f(m)\otimes_A c) = f(m)\sw 0\cdot\gamma(f(m)\sw 1\otimes_A c) = 
f(m\sw
0)\cdot\gamma(m\sw 1\otimes_A c) = f\circ\nu_M(m\otimes_A c),
$$
where we used the fact that $f$ is a right $\cC$-comodule map. Finally
we take any $m\in M$ and compute
$$
\nu_M\circ \Phi_M(m) = \nu_M(m\sw 0\otimes_A m\sw 1) = m\sw
0\cdot\gamma(m\sw 1\otimes_A m\sw 2) = m\sw 0\cdot \eps_\cC(m\sw 1) = m.
$$
This shows that $\nu$ is the required splitting of $\Phi$.
\end{proof}

Recall from \cite{Guz:coi} that an $A$-coring $\cC$ is said to be 
{\em coseparable} if there exists a $(\cC,\cC)$-bicomodule splitting of
the coproduct. Explicitly one requires an $(A,A)$-bimodule map 
$\pi :\cC\otimes _A\cC\to \cC$
such that 
\begin{equation}
(\cC\otimes_A\pi)\circ(\Delta_\cC\otimes_A\cC) = 
\Delta_\cC\circ\pi = (\pi\otimes_A\cC)\circ(\cC\otimes_A\Delta_\cC),
\quad \pi\circ\Delta_\cC = \cC.
\label{cosep}
\end{equation}
Theorem~\ref{thm.sep2} implies the following characterisation of
coseparable corings which supplements
\cite[Theorem~3.10]{Guz:coi}. 
\begin{corollary}
An $A$-coring $\cC$ is coseparable if and only if the forgetful functor 
$F:\M^\cC\to\M_A$ is separable.
\label{cor2}
\end{corollary}
\begin{proof}
Given a map $\gamma$ as in Theorem~\ref{thm.sep2} one defines
$\pi:\cC\otimes_A\cC \to \cC$, $c\otimes_A c'\mapsto c\sw 1\cdot\gamma(c\sw
2\otimes_Ac')$. By the ``$\Leftarrow$" part of the proof of
Theorem~\ref{thm.sep2}, $\pi = \nu_\cC$ and thus it 
is a right $\cC$-comodule splitting of
$\Delta_\cC$. Using the fact that $\Delta_\cC$ is an $(A,A)$-bimodule
map one easily verifies that $\pi$ is a left $\cC$-comodule map.
Conversely, given $\pi$ define $\gamma = \eps_\cC\circ\pi$. Since $\pi$
splits $\Delta_\cC$, $\gamma\circ\Delta_\cC = \eps_\cC$. Applying
$\cC\otimes_A \eps_\cC$ to the first equality in (\ref{cosep}) and
$\eps_\cC\otimes_A\cC$ to the second equality in (\ref{cosep}) one
deduces eq.\ (\ref{eq}). By Theorem~\ref{thm.sep2}, the forgetful
functor is separable as required.
\end{proof}

\begin{corollary}
Let $B\to A$ be a ring extension such that $A$ is faithfully
flat as either left or right $B$-module, and let $\cC$ 
 be the canonical
$A$-coring in Example~\ref{can.coring}. The extension $B\to A$ is split if and only if the
functor $F:\M^\cC\to\M_A$ is separable.
\end{corollary}
\begin{proof}
Recall that an extension $B\to A$ is split iff there exists a
$(B,B)$-bimodule map $E:A\to B$ such that $E(1)=1$ (cf.\
\cite{Pie:ass}).
 In the case of the
canonical $A$-coring $A\otimes_BA$ the conditions required for the map
$\gamma \in {}_A{\rm Hom}_A (A\otimes_BA\otimes_BA, A)$ read 
$\gamma(a\otimes_B 1\otimes_B a')= aa'$ and $a\otimes_B
\gamma(1\otimes_B a'\otimes_Ba'') = \gamma(a\otimes_Ba'\otimes_B
1)\otimes_B a''$ for all $a,a',a''\in A$. Since 
${}_A{\rm Hom}_A (A\otimes_BA\otimes_BA, A)
\cong {}_B{\rm Hom}_B (A, A)$ the maps $\gamma$ are in one-to-one
correspondence with the maps $E\in{}_B{\rm Hom}_B (A, A)$ via
$\gamma(a\otimes_Ba'\otimes_Ba'') = aE(a')a''$. The first of the above
conditions for $\gamma$ is equivalent to the normalisation of $E$, $E(1)
=1$, while the second condition gives for all $a\in A$, $1\otimes_B E(a)
= E(a)\otimes_B 1$. By the faithfully flat descent the latter is
equivalent to $E(a)\in B$.
\end{proof}

In the case of the coring $\cC = A\otimes C$ corresponding to an entwining
structure $(A,C)_\psi$ via Proposition~\ref{prop.tak}, a map $\gamma$ in
Theorem~\ref{thm.sep2} corresponds to a normalised integral map  
in $(A,C)_\psi$ in
the sense of \cite[Definition~4.1]{Brz:coa}. Thus Theorem~\ref{thm.sep2}
implies \cite[Theorem~4.2]{Brz:coa}, the latter being a generalisation of
\cite[Theorem~2.3]{CaeMil:sep} formulated for Doi-Koppinen modules.

As explained in \cite{CaeMil:sep} (cf.\ \cite[Proposition~1.2]{Nas:sep})
 the separability of the forgetful
functor implies various Maschke-type theorems. Thus, we have the
following generalisation of \cite[Corollary~4.4]{Brz:coa}

\begin{corollary}{\rm\bf (Maschke-type theorem for a coring)}
Let $\cC$ be an $A$-coring with a map $\gamma$ satisfying hypothesis of
Theorem~\ref{thm.sep2}. Then every right $\cC$-comodule
which is semisimple (resp.\ projective, injective) as a right $A$-module  is
semisimple (resp.\ projective, injective) $\cC$-comodule.
\end{corollary}

\section{Frobenius properties of a coring}
Let $\cC$ be an $A$-coring. As explained in
\cite[3.2~Proposition~(a)]{Swe:pre}, $R = {}_A{\rm Hom}(\cC,A)$ is a
ring with unit $\eps_\cC$ and product $(rr')(c) = r'(c\sw 1\cdot r(c\sw
2))$, for all $r,r'\in R$ and $c\in \cC$. $R$ is a left $A$-module via
$(a\cdot r)(c) = r(c\cdot a)$, for all $a\in A, c\in\cC, r\in R$. 
Furthermore the map $\iota :
A\to R$ given by $\iota(a)(c) = \eps_\cC(c)a$ is a ring map. In this
section we study when $A\to R$ is a Frobenius extension.

Recall from \cite{Kas:pro} and \cite{NakTsu:fro} that a ring
 extension $A\to R$
is called a {\em Frobenius extension} (of the first kind) iff $R$ is a
finitely generated projective right $A$-module and $R \cong 
{\rm Hom}_A(R,A)$ as $(A,R)$-bimodules. The $(A,R)$-bimodule
structure of ${\rm Hom}_A(R,A)$ is given by $(a\cdot f\cdot r)(r') =
af(rr')$, for all $a\in A, r,r'\in R$ and $f\in {\rm Hom}_A(R,A)$.
Frobenius extensions are closely related with a certain type of corings,
namely, $A\to R$ is a Frobenius extension if and only if $R$ is an
$A$-coring such that the coproduct is an $(R,R)$-bimodule map (cf.\
\cite[Proposition~4.3]{Kad:new}, \cite[Remark~2.5]{CaeIon:fro}). The
main result of this section is contained in the following
\begin{theorem}
Let $\cC$ be an $A$-coring and let $R = {}_A{\rm Hom}(\cC,A)$. If $\cC$
is a projective left $A$-module then the following are equivalent:

(1) The forgetful functor $F:\M^\cC\to \M_A$ is a Frobenius functor. 

(2) $\cC$ is a finitely generated left $A$-module and 
the ring extension $A\to R$ is Frobenius.

(3) $\cC$ is a finitely generated left $A$-module and $\cC\cong R$ as
$(A,R)$-bimodules, where $\cC$ is a right $R$-module via $c\cdot r =
c\sw 1\cdot r(c\sw 2)$, $\forall c\in \cC, r\in R$.

(4) $\cC$ is a finitely generated left $A$-module and  there exists
$e\in \cC^A$ such that the map $\phi: R\to \cC$, $r\mapsto e\sw 1\cdot
r(e\sw 2)$ is bijective.
\label{thm.fro}
\end{theorem}

Recall from \cite{CaeMil:Doi} that a functor is said to be {\em Frobenius} in case it
 has the same 
right and left adjoint. Since by Lemma~\ref{lemma.adjoint} 
the functor $-\otimes_A\cC:\M_A\to \M^\cC$ is the right adjoint of
the forgetful functor $F$,
Theorem~\ref{thm.fro}~(1) is equivalent to the statement that 
 $-\otimes_A\cC$ is the left adjoint of
 $F$. 
Proof of  Theorem~\ref{thm.fro} is based on two lemmas 
\begin{lemma}
If  
$-\otimes_A\cC:\M_A\to \M^\cC$ is the left adjoint of
the forgetful functor $F:\M^\cC\to \M_A$, then $\cC$
is a finitely generated left $A$-module. 
\label{lemma.fro}
\end{lemma}
\begin{proof} The proof of this lemma is based on the proofs of
\cite[Lemma~2.3, Theorem~2.4 1) $\Rightarrow$ 2)]{CaeMil:Doi}. Let for
all $M\in \M_A$, $N\in\M^\cC$, $\eta_{M,N} :{\rm Hom}^\cC(M\otimes_A\cC,
N) \to {\rm Hom}_A(M,N)$ be the natural isomorphism and let $e =
\eta_{M,N}(\cC)(1)$. Since $\eta$ is natural and $A\in \M_A$, $\cC\in
\M^\cC$ we have for all $f\in {\rm Hom}^\cC(\cC,
\cC)$, $\eta_{A,\cC}\circ {\rm Hom}^\cC(\cC, f) = {\rm Hom}_A(A, f)\circ
\eta_{A,\cC}$. Evaluating this equality at $\cC$ and the resulting
equality  at $1$  we obtain $\eta_{A,\cC}(f)(1) = f(\eta_{A,\cC}(\cC)(1)) =
f(e)$. Now, for any $c\in\cC$, let $f_c\in {\rm Hom}^\cC(\cC,
\cC)$ be the unique morphism in $\M^\cC$ such that $\eta_{A,\cC}(f_c)(a)
= c\cdot a$, $\forall a\in A$. Taking $a = 1$ and using above equality we
obtain $c= \eta_{A,\cC}(f_c)(1) = f_c(e)$. Let $\Delta_\cC(e) =
\sum_{i=1}^n 
e_i\otimes_A \bar{e}_i$. Since $f_c$ is a right $\cC$-comodule map we
have
$
c\sw 1\otimes_A c\sw 2 = \sum_{i=1}^n f_c(e_i)\otimes_A\bar{e}_i$. Applying
$\eps_\cC\otimes_A \cC$ to this equality we obtain for all $c\in \cC$, $c
= \sum_{i=1}^n\eps_\cC(f_c(e_i))\cdot \bar{e}_i$, i.e., $\cC$ is a finitely
generated left $A$-module.
\end{proof}

In fact using the same techniques as in the proof of 
\cite[Lemma~2.3]{CaeMil:Doi} one can establish that with the notation
and hypothesis of Lemma~\ref{lemma.fro}, $\eta_{M,N}(f)(m) =
f(m\otimes_A e)$, for all $f\in {\rm Hom}^\cC(M\otimes_A\cC, N)$, $m\in
M$. 

\begin{lemma}
Let $\cC$ be an $A$-coring and let $R = {}_A{\rm Hom}(\cC,A)$. If $\cC$
is a finitely generated projective  left $A$-module then the categories
$\M^\cC$ and $\M_R$ are isomorphic with each other.
\label{lemma.iso}
\end{lemma}
\begin{proof} 
Given $M\in \M^\cC$ one can view it as a right $R$-module via $m\cdot r
= m\sw 0\cdot r(m\sw 1)$, for all $m \in M$, $r\in R$. Indeed, it is
immediate that $m\cdot\eps_\cC = m$. Furthermore, using that the right
coaction of $\cC$ on $M$ is a right $A$-module map we have for all
$r,r'\in R$, $m\in M$,
$$
(m\cdot r)\cdot r' = (m\cdot r)\sw 0\cdot r'((m\cdot r)\sw 1) = m\sw
0\cdot r'(m\sw 1\cdot r(m\sw 2)) = m\sw 0\cdot (rr')(m\sw 1) = m\cdot
rr'.
$$
Let $\{r_i, c^i\}_{i=1}^n$, $r_i\in R$, $c^i\in \cC$ 
be a dual basis of $\cC$ as a left $A$-module.
Notice that for all $c\in \cC$, 
$\Delta_\cC (c) = \sum_ir_i(c)\cdot c^i\sw 1\otimes_A c^i\sw 2 $. On the
other hand
\begin{eqnarray*}
\Delta_\cC(c) &=& \sum_ic\sw 1\otimes_A r_i(c\sw 2)\cdot c^i = 
\sum_ic\sw 1\cdot r_i(c\sw 2)\otimes_A c^i \\
&=& \sum_{i,j} r_j(c\sw 1\cdot r_i(c\sw 2))\cdot c^j\otimes_A c^i = 
\sum_{i,j}r_ir_j(c)\cdot
c^i\otimes_A c^j.
\end{eqnarray*}
This implies that
$$
\sum_ir_i\otimes_Ac^i\sw 1\otimes_A c^i\sw 2 = \sum_{i,j}r_ir_j\otimes
_Ac^i\otimes_A c^j .
$$ 
Using this equality one easily finds that given $M\in\M_R$, $M$ is a
right $\cC$-comodule with the coaction $\rho^M(m) = \sum_im\cdot
r_i\otimes_A c^i$. Also, one easily checks that the maps described
provide the required isomorphism of categories.
\end{proof}

Now we can  prove Theorem~\ref{thm.fro}

\begin{proof}
(1) $\Leftrightarrow$ (2). By Lemma~\ref{lemma.fro}, (1) implies that 
$\cC$ is a finitely generated projective left $A$-module. 
Lemma~\ref{lemma.iso} shows that the forgetful functor is the restriction
of scalars functor $\M_R\to \M_A$. By Lemma~\ref{lemma.adjoint}, this
functor has the right adjoint $-\otimes_A\cC$. 
By \cite[Theorem~3.15]{MenNas:ind} the restriction of scalars
functor has the same left and right adjoint if and only if the
extension $A\to R$ is Frobenius.

(2) $\Leftrightarrow$ (3). Since $\cC$ is a finitely generated
projective left $A$-module the map $\alpha :\cC\to{\rm Hom}_A(R,A)$,
given by $\alpha(c)(r) = r(c)$ is bijective by \cite[3.5~Duality
Lemma]{Swe:pre}. Clearly $\alpha$ is an $(A,R)$-bimodule map. Thus we
have $\cC \cong {\rm Hom}_A(R,A)$ as $(A,R)$-bimodules. The extension $A\to R$ is Frobenius
iff $R\cong {\rm Hom}_A(R,A)$, i.e., iff $\cC \cong R$ as
$(A,R)$-bimodules. 

(3) $\Leftrightarrow$ (4). This follows from the bijective
correspondence $\theta: \cC^A\to {}_A{\rm
Hom}_R(R, \cC)$, $\theta(e)(r) = e\cdot r = e\sw 1\cdot r(e\sw 2)$,
$\theta^{-1}(f) = f(\eps_\cC)$.
%
\end{proof}

In the case of the canonical coring associated to a ring extension $B\to
A$, Theorem~\ref{thm.fro} gives the criteria when the extension $A\to
{}_B{\rm End}(A)$ is Frobenius. The ring structure on ${}_B{\rm End}(A)$
is given by the opposite composition of maps. 
For example if $B\to A$ is
itself a Frobenius extension with a Frobenius system $E\in {}_B{\rm
Hom}_B(A,B)$, $a_i,\overline{a}^i\in A$, $i=1,\ldots, n$ (i.e., for all
$a\in A$, $\sum_ia_iE(\overline{a}^ia) = \sum_iE(aa_i)\overline{a}^i
=a$), then  $A\to R$ is Frobenius
by the endomorphism ring theorem \cite{Kas:pro}. In this case  
$e=\sum_ia_i\otimes_B\overline{a}^i$ and the inverse of $\phi$ is given
by $\phi^{-1}(a\otimes_Ba')(a'') = E(a''a)a'$, for all $a,a',a''\in A$. 
In the case of the coring associated to 
an entwining structure
$(A,C)_\psi$, $R$ becomes a generalised smashed product
$C^{*op}\#_{\bar{\psi}} A$, and therefore Theorem~\ref{thm.fro} implies
\cite[Theorem~2.4]{CaeMil:Doi} or its entwined module formulation 
\cite[Proposition~3.5]{Brz:fro}.

\section{Galois-type corings}
In this section we study $A$-corings $\cC$ for which $A\in\M^\cC$.
\begin{lemma}
For an $A$-coring $\cC$, $A$ is a right $\cC$ comodule if and only if
there exists a grouplike $g\in\cC$ (cf.\ \cite[1.7~Definition]{Swe:pre}).
\label{lemma.grouplike}
\end{lemma}
\begin{proof}
If $A\in\M^\cC$ define $g=\rho^A(1) \in A\otimes_A\cC \cong \cC$. Then
$$
\Delta_\cC(g) = \Delta_\cC(\rho^A(1)) =
(A\otimes_A\Delta_\cC)\circ\rho^A(1) =
(\rho^A\otimes_A\cC)\circ\rho^A(1) = \rho^A(1)\otimes_A\rho^A(1) =
g\otimes_Ag .
$$
Furthermore $\eps_\cC(g) = (A\otimes_A\eps_\cC)\circ\rho^A(1) =1$, so
that $g$ is a grouplike as required.

Conversely, if $g\in\cC$ is a grouplike, define $\rho^A:A\to \cC\cong
A\otimes_A\cC$, $a\mapsto g\cdot a = 1\otimes_Ag\cdot a$. Clearly
$\rho^A$ is a right $A$-module map. The fact that $\eps_\cC$ is a right
$A$-module map implies for all $a\in A$,
$(A\otimes_A\eps_\cC)\circ\rho^A(a) = \eps_\cC(g)a = a$. Finally, since
$\Delta_\cC$ is a right $A$-module map we have
$(A\otimes_A\Delta_\cC)\circ\rho^A(a) = 1\otimes_A \Delta_\cC(g)\cdot a
=1\otimes_A 
g\otimes_Ag\cdot a$. On the other hand $(\rho^A\otimes_A\cC)\circ\rho^A(a) = 
\rho^A(1)\otimes_A g\cdot a = 1\otimes_Ag\otimes_Ag\cdot a$. Put together
this implies that $\rho^A$ is a right $\cC$-coaction.
\end{proof}

For example, if $\cC$ is the canonical $A$-coring associated to a ring 
extension $B\to A$ then $g=1\otimes_B1$ is a grouplike by
\cite[1.9~Proposition~(a)]{Swe:pre}, and hence $A$ is a right $\cC$-comodule
via $a\mapsto 1\otimes_B a$. 

For the rest of this section we assume that
$A \in \M^\cC$ and denote the corresponding grouplike by $g$. In this
case, for each $M\in\M^\cC$ we define the {\em coinvariants} by
$$
M^{co\cC} = \{m\in M\; |\; \rho^M(m) = m\otimes_A g\}.
$$
In particular, let  
$B = A^{co\cC} = \{b\in A\;|\; b\cdot g = g\cdot b\}$, i.e., $B$ is the
centraliser of $g$ in $A$. Clearly $B$
is a subring of $A$, $\rho^A$ is a $(B,A)$-bimodule map and
$M^{co\cC}\in \M_B$.
\begin{proposition}
Let $\cC$ be an $A$-coring with a grouplike $g$, and let $B=A^{co\cC}$. 
Then the functor $G:\M^\cC\to \M_B$, $M\mapsto M^{co\cC}$ is the right
adjoint of the induction functor $-\otimes_BA:\M_B\to\M^\cC$. Here, for any $M\in
\M_B$, the right $\cC$-coaction on $M\otimes_BA$ is given by
$M\otimes_B\rho^A$. 
\label{adjoint2}
\end{proposition}
\begin{proof}
For any $M\in\M^\cC$ define an additive map 
$\Psi_M: M^{co\cC}\otimes_B A\to M$,
$m\otimes_B a\mapsto m\cdot a$.  Clearly
$\Psi_M$ is a right $A$-module map. Furthermore for any $m\in M^{co\cC}$
and $a\in A$ we have
$\rho^M(m\cdot a) = m\otimes_Ag\cdot a$ as well as
$$
(\Psi_M\otimes_A\cC)\circ\rho^{G(M)\otimes_BA} (m\otimes_B a) =
(\Psi_M\otimes_A\cC)(m\otimes_B1\otimes_Ag\cdot a) = m\otimes_Ag\cdot a.
$$
Therefore $\Psi_M$ is a $\cC$-comodule map. One easily checks that
$\Psi_M$ is natural in $M$. 

Next for any $N\in\M_B$ define
an additive map $\Phi_N:N\to (N\otimes_BA)^{co\cC}$, 
$n\mapsto n\otimes_B 1$. Notice that $(N\otimes_B \rho^A)(n\otimes_B 1)
= n\otimes_B 1\otimes_A g$ so that $\Phi_N$ is well-defined. 
One easily checks that $\Phi_N$ is a right $B$-module map natural in
$N$. 

Finally, with $M$, $N$ as before, take any $m\in M^{co\cC}$ and compute
$\Psi_M\circ\Phi_{M^{co\cC}}(m) = \Psi_M(m\otimes_B 1) = m$. Then take
any $n\in N$, $a\in A$ and compute $\Psi_{N\otimes_B
A}\circ(\Phi_N\otimes_BA)(n\otimes_B a) = \Psi_{N\otimes_B
A}(n\otimes_B1\otimes_B a) = n\otimes_Ba$. This proves that $\Phi$ and
$\Psi$ are the required adjunctions.
\end{proof}

In fact Proposition~\ref{adjoint2} is a special case of the following
Hom-Tensor relation for a coring. Let $B,A$ be rings, $\cC$ an
$A$-coring , and let $V$ be a
$(B,A)$-bimodule and right $\cC$-comodule such that the coaction is a
left $B$-module map. Then for any right $B$-module $N$ and a right
$\cC$-comodule $M$ one has 
$$
{\rm Hom}^\cC(N\otimes_B V, M)\cong {\rm Hom}_B(N, {\rm Hom}^\cC(V,M)).
$$
Explicitly for any $f\in {\rm Hom}^\cC(N\otimes_B V, M)$ , $f\mapsto
(n\mapsto (v\mapsto  f(n\otimes_Bv)))$ and for any $g\in {\rm Hom}_B(N, {\rm
Hom}^\cC(V,M))$, $g\mapsto (n\otimes_B v\mapsto g(n)(v))$.

\begin{definition}
Let $\cC$ be an $A$-coring with a grouplike $g$, and let $B=A^{co\cC}$.
$\cC$ is said to be {\em Galois} iff there exists an $A$-coring
isomorphism $\chi: A\otimes_BA\to \cC$ such that $\chi(1\otimes_B1) =g$.
\end{definition}
For example, if $A$ is a division ring, any $A$-coring  which is 
generated by a grouplike $g$ as an $(A,A)$-bimodule is Galois 
(cf.\ \cite[2.2~Fundamental Lemma]{Swe:pre}).
Our present terminology is motivated by the following two examples.
\begin{example}
Let $(A,C)_\psi$ be an entwining structure and let $\cC = A\otimes C$ as
in Proposition~\ref{prop.tak}. $\cC$ is a Galois $A$-coring if and only
if $(A,C)_\psi$ is the canonical entwining structure of a $C$-Galois
extension $B\subset A$.
\label{ex.galois}
\end{example}
\begin{proof}
 If $\cC$ corresponds to a $C$-Galois extension $B\subset A$ then
$A\otimes_B A\cong A\otimes C$ as $A$-corings via the canonical map
$can$, and hence $\cC$ is Galois. Conversely, if $\cC$ is Galois then
$A$ is a right $\cC$-comodule, and by the correspondence in
Proposition~\ref{prop.tak} it is an $(A,C)_\psi$-module. The
corresponding grouplike in $\cC$ is $g=\rho^A(1)= 1\sw 0\otimes 1\sw 1$.
Furthermore $A^{co\cC} = \{b\in A \; | \rho^A(b) = b1\sw 0\otimes 1\sw
1\} = A^{coC}$, since $A\in\M_A^C(\psi)$. For the same reason the 
$A$-coring
isomorphism $\chi :A\otimes_BA\to A\otimes C$ explicitly reads
$\chi(a\otimes_Ba') = a\cdot(1\sw 0\otimes 1\sw 1)\cdot a' = a1\sw 0a'_\alpha\otimes 1\sw 1^\alpha =
aa'\sw
0\otimes a'\sw 1$ and thus coincides with the canonical map $can$. This
proves that $B\subset A$ is a $C$-Galois extension and by the uniqueness
of the canonical entwining structure, $(A,C)_\psi$ must be the canonical
entwining structure associated to $B\subset A$.
\end{proof}
\begin{example}
Let $(A,C,\psi)$ be a weak entwining structure corresponding to a weak
$C$-Galois extension $B\subset A$ as described in Example~\ref{ex.weak}.
Then the corresponding $A$-coring $\cC \subset A\otimes C$ given in
Proposition~\ref{weak} is Galois.
\label{weak.galois}
\end{example}
\begin{proof}
It suffices to show that ${\rm Im}(can) = \cC$, then $can$ will provide
the required isomorphism of $A$-corings. Notice that from the definition
of $\psi$ in Example~\ref{ex.weak} 
it follows that ${\rm Im}\psi \subseteq {\rm Im}(can)$. Since
a typical element of $\cC$ is of the form 
$a1_\alpha\otimes c^\alpha$ and $can$ is a left
$A$-module map we have $a1_\alpha\otimes c^\alpha\in {\rm Im}(can)$.
Therefore $\cC\subseteq {\rm Im} (can)$. On the other hand, since $A$ is
a weak entwined module we have for all
$a\in A$, $\rho^A(a) = \rho^A(a1) = 
a\sw 01_\alpha\otimes a\sw 1^\alpha\in \cC$.
In the view of the fact that $can(a\otimes_Ba') = a\rho^A(a')$ this
implies that ${\rm Im}(can)\subseteq \cC$.
\end{proof}
\begin{theorem}
Let $\cC$ be an $A$-coring with a grouplike $g$, $B=A^{co\cC}$, and let 
$G:\M^\cC\to \M_B$, $M\mapsto M^{co\cC}$. 
If $\cC$ is Galois and $A$ is a faithfully flat left $B$-module then 
the functors $G$ and $-\otimes_BA:\M_B\to \M^\cC$ 
are inverse equivalences. Conversely, if $G$ and $-\otimes_BA$ 
are inverse equivalences then $\cC$ is Galois. In this case if 
$\cC$ is a flat left $A$-module then $A$ is a faithfully flat left
$B$-module.
\label{thm.galois}
\end{theorem}
\begin{proof}
Assume that $\cC$ is Galois and $A$ is a faithfully flat left
$B$-module. First notice that $\chi: a\otimes_Ba'\mapsto a\cdot
g\cdot a'$. For all $M\in\M^\cC$ consider the following commuting
diagram of right $\cC$-comodule maps
$$
\begin{CD}  
0 @>>>   M^{co\cC}\otimes_B A  @>>>  
M\otimes_B A 
  @>>>( M\otimes_A\cC)\otimes_B{A} \\  
  @.      @VV\Psi_{M}V  
@VV{M\otimes_{A}\chi}V 
@VV{({M}\otimes_A{\cC})\otimes_{{A}}\chi}V\\  
0 @>>>  {M}   @>\rho^{M}>>   {M} \otimes_A {\cC}  
  @>\ell_{{M}{\cC}}>> {M}\tens_A {\cC}\otimes_A {\cC}  
\end{CD}    
$$
The maps in the top row are the obvious inclusion and $m\otimes_Ba\mapsto
\rho^M(m)\otimes_B a - m\otimes_A g\otimes_B a$, while 
$\ell_{M\cC}=\rho^M\otimes_A\cC - M\otimes_A\Delta_\cC$ 
is the coaction equalising map. The top row is exact since it is a
defining sequence of $M^{co\cC}$ tensored with $A$, and $-\otimes_BA$ is
exact. The bottom row is exact too. Since $\chi$ is a bijection, so are
$M\otimes_A\chi$ and $(M\otimes_A\cC)\otimes_A\chi$. Therefore $\Psi_M$
is an isomorphism in $\M^\cC$.

For all  $N\in\M_B$ consider the following commutative diagram of right
$B$-module maps
$$
\begin{CD}
0@>>> N @>>> N\otimes_B A @>>> N\otimes_B A\otimes_B A\\
@. @VV{\Phi_N}V @VV{=}V @VV{N\otimes_B\chi}V\\
0@>>> (N\otimes _BA)^{co\cC} @>>> N\otimes_B A @>>>
N\otimes_BA\otimes_A\cC
\end{CD}
$$
The maps in the top row are: $n\mapsto n\otimes_B 1$ and $n\otimes_B
a\mapsto n\otimes_B a\otimes_B 1-n\otimes_B1\otimes_B a$, and the top 
row is exact by the faithfully flat descent. The bottom row is the 
defining sequence of
$(N\otimes_B A)^{co\cC}$ and hence is exact. This implies that 
$\Phi_N$ is an isomorphism in $\M_B$ and
completes the proof of the fact that $G$ and $-\otimes_BA$ are inverse
equivalences. 

Conversely, assume that $G$ and $-\otimes_BA$ are inverse equivalences. 
Notice that $\phi:A\to\cC^{co\cC}$, $a\mapsto a\cdot g$ is an
$(A,B)$-bimodule isomorphism. Since $\cC$ is a right $\cC$-comodule via
the coproduct, there is a corresponding adjunction
$\Psi_{\cC}:\cC^{co\cC}\otimes_BA\to \cC$, and it is bijective. Define
$\chi:A\otimes_BA\to \cC$, $\chi = (\phi\otimes_BA)\circ\Psi_\cC$.
Explicitly $\chi:a\otimes_Ba'\mapsto a\cdot g\cdot a'$. Clearly $\chi$
is an $(A,A)$-bimodule map such that $\chi(1\otimes_B1) = g$.
Furthermore since $g$ is a  grouplike and $\Delta_\cC$ is an
$(A,A)$-bimodule map we have $\Delta_\cC\circ\chi(a\otimes_B a') = a\cdot
g\otimes_Ag\cdot a'$, for all $a,a'\in A$. On the other hand
$$
(\chi\otimes_A\chi)\circ\Delta_{A\otimes_B A}(a\otimes_Ba') = \chi(a\otimes_B
1)\otimes_A\chi(1\otimes_B a') = a\cdot g\otimes_A g\cdot a'.
$$
Therefore $\chi$ is an $A$-coring isomorphism and hence $\cC$ is Galois.

If $\cC$ is a flat left $A$-module then both kernels and cokernels of
any morphism in $\M^\cC$ are right $\cC$-comodules (i.e., $\M^\cC$ is an
Abelian category). Therefore any sequence of $\cC$-comodule maps is
exact if and only if it is exact as a sequence of additive maps. Since 
 $-\otimes_BA$ is an equivalence, it preserves and reflects exact
sequences. By the above observation it does so even as viewed as a
functor from $\M_B$ to the category of $\bf Z$-modules. 
Therefore $A$ is a faithfully flat left $B$-module.
\end{proof}

In view of Example~\ref{ex.galois}, Theorem~\ref{thm.galois} is a 
generalisation of
\cite[Corollary~3.11]{Brz:mod} which in turn is a generalisation of
\cite[Theorem~3.7]{Sch:pri} (from which the idea of the proof is taken).

\section{Comments on duality and outlook}
Throughout this section we assume that $k$ is a field. The results of the
previous sections have dual counterparts. To introduce them we first
propose the following dualisation of the notion of a coring. 

Let $C$ be a coalgebra. A {\em $C$-algebroid}  (or a {\em $C$-ring}) 
is a $(C,C)$-bicomodule
$\cA$ together with $(C,C)$-bicomodule maps $\mu_\cA:\cA\square_C\cA \to
\cA$ and $\eta_\cA:C\to \cA$ such that 
$$
\mu_\cA\circ(\mu_\cA\square_C\cA) = \mu_\cA\circ(\cA\square_C \mu_\cA),
\quad \mu_\cA\circ(\eta_\cA\square_C\cA)\circ{}^\cA\rho = 
\mu_\cA\circ(\cA\square_C \eta_\cA)\circ\rho^\cA = \cA.
$$
Here $\square_C$ denotes the cotensor product of $(C,C)$-bicomodules,
and ${}^\cA\rho$, $\rho^\cA$ are left, right coactions of $C$ on $\cA$. 
The map $\mu_\cA$ is called a product and $\nu_\cA$ is called a unit of
the $C$-algebroid $\cA$. Obviously a $k$-algebroid is simply a
$k$-algebra. An example of a $C$-algebroid can be obtained by
dualisation of Example~\ref{can.coring}.
\begin{example}
Let $\pi:C\to B$ be a morphism of coalgebras, then $C$ is a
$(B,B)$-comodule via $(C\otimes \pi)\circ\Delta$ and $(\pi\otimes
C)\circ\Delta$. Define $\cA = C\square_BC$. Then $\cA$ is a
$C$-algebroid with product $\mu_\cA :A\square_C\cA \cong
C\square_BC\square_BC\to \cA$, $\mu_\cA = C\square_B\eps\square_B C$ and
unit $\nu_A = \Delta$.
\end{example}
A right $\cA$ module is a right $C$-comodule $M$ with a
map $\rho_M :M\square_C\cA\to M$ such that
$$
\rho_M\circ(\rho_M\square_C\cA) = \rho_M\circ(M\square_C \mu_\cA),
\quad  
\rho_M\circ(M\square_C \eta_\cA)\circ\rho^M = M,
$$
where $\rho^M$ is the right coaction of $C$ on $M$. 
The map $\rho_M$ is called the action of $A$ on $M$. A map of $\cA$-modules is a right
$C$-colinear map which respects the actions. The category of right 
$\cA$-modules is denoted by $\M_\cA$. 

Dualising Proposition~\ref{prop.tak} one obtains 
\begin{proposition}
For an entwining structure $(A,C)_\psi$, view $\cA = C\otimes A$ as a
$(C,C)$-bicomodule with the left coaction ${}^\cA\rho = \Delta\otimes C$
and the right coaction 
$\rho ^\cA = (C\otimes \psi)\circ (\Delta\otimes A)$.  Then 
$\cA$ is a
$C$-algebroid   with the product
 $\mu_\cA : \cA\square_C\cA \cong C\otimes A\otimes A\to \cA$, 
$\mu_\cA = C\otimes \mu$,  and the counit 
$\eta_\cA = C\otimes 1$.

Conversely if $C\otimes A$ is a $C$-algebroid 
with the product and the
unit above and the natural left $C$-comodule structure $\Delta\otimes A$, 
then $(A,C)_\psi$ is an entwining structure, where 
$\psi = (\eps\otimes A\otimes C)\circ \rho^{C\otimes A}$.

Under this bijective correspondence $\M_\cA = \M_A^C(\psi)$.
\label{prop.tak*}
\end{proposition}
One can also dualise main results in Section~3. For example, given a
right $C$-comodule $M$ one defines a right $\cA$-module structure on
$M\square_C\cA$ via $M\square_C\mu_\cA$. Dualising
Lemma~\ref{lemma.adjoint} one easily finds that the functor
$-\square_C\cA: \M^C\to \M_\cA$ is left adjoint of the forgetful functor
$\M_\cA\to \M^C$. One can then study  separability of these functors
and thus obtain results dual to Theorem~\ref{thm.sep1} and
Theorem~\ref{thm.sep2}: 
\begin{theorem}
The functor $-\square_C\cA: \M^C\to \M_\cA$ is separable if and only if
there exists a map $e: \cA\to k$ such that $(e\otimes
C)\circ\rho^\cA = (C\otimes e)\circ{}^\cA\rho$ and $e\circ\eta_\cA
= \eps$.

The forgetful functor $\M_\cA\to \M^C$ is separable if and only if there
exists $\gamma\in{}^C{\rm Hom}^C(C, \cA\square_C\cA)$ such that
$$
(\mu_\cA\square_C\cA)\circ(\cA\square_C\gamma)\circ\rho^\cA =
(\cA\square_C\mu_\cA)\circ(\gamma\square_C\cA)\circ{}^\cA\rho.
$$
and $\mu_\cA\circ\gamma =\eta_\cA$.
\end{theorem}
Finally one can proceed to formulate results dual to ones described in
Sections~4 and 5. In particular one easily finds that $C$ is an $\cA$-module
provided there is a nontrivial character $\kappa:\cA\to k$. Using the 
natural identification
$C\square_C \cA \cong \cA$ the action  of $\cA$ on $C$ is a map $\rho_C
:\cA\to C$, $a\mapsto \kappa(a\sw 0)a\sw 1$. In this case
one can study invariants of $C$, $I = span \{\kappa(a\sw 0)a\sw 1 -
a\sw{-1}\kappa(a\sw 0) \; |\; a\in \cA\}$, where ${}^\cA\rho(a) =
a\sw {-1}\otimes a\sw 0$ denotes the left coaction of $C$ on $\cA$, 
and find that $I$ is a coideal. Hence $B = C/I$
is a coalgebra. One can proceed to define a {\em Galois $C$-algebroid}
as a $C$-algebroid $\cA$ with a nontrivial character $\kappa$ such that
there is an isomorphism of $C$-algebroids $\chi :\cA\to C\square_B C$
such that $\kappa = (\eps\square_B\eps)\circ\chi$. The details of all the
results quoted here as well as derivation of results 
concerning functor $-\square_B \cA :\M^B\to \M_\cA$ are left to
the reader.

In this paper we have shown that corings provide 
a natural framework for studying (weak) entwining structures and
(weak) entwined
modules, and that some general results for the latter can be deduced
from the corresponding results for corings. It seems natural and,
indeed, desired to study other properties of entwining structures from
the point of view presented in this paper. On the other hand already
known properties of entwining structures might suggest similar
properties for corings. For example, we expect that the pairs of adjoint
functors considered in this paper are in fact special cases of
 pairs of induction functors that involve both a coring and
an algebroid. Another line of development is suggested by the following
observation made by S.\ Caenepeel
\cite{Cae:pri}. The relationship between corings and weak entwining 
structures described
in Proposition~\ref{weak} can be made more complete along the lines
similar to Proposition~\ref{prop.tak} provided one introduces the
following generalisation of a coring. For a
ring $A$, an $A$-{\em pre-coring} is an $(A,A)$-bimodule $\cC$
non-unital as a right $A$-module (i.e., it is not required that $c\cdot 1
= c$ for all $c\in \cC$) with $(A,A)$-bimodule maps  
$\Delta_\cC:\cC\to \cC\otimes_A\cC$ and 
 $\eps_\cC:\cC\to A$ satisfying the
same axioms as for a coring and such that $c\cdot 1= 
\eps_\cC(c_{(1)}\cdot 1)c_{(2)}$  for all  $c\in \cC$.
S.\ Caenepeel then proves that given an $A$-pre-coring $\cC$, 
the map $p:\ {\cal C}\to {\cal C}$
given by $p(c)=c\cdot 1$ is an $(A,A)$-bilinear projection, 
$p\circ p=p$ such that $(p\otimes_A p)\circ \Delta_\cC= \Delta_\cC \circ
p$. Furthermore
$\overline{\cal C}={\rm Coim}(p)\cong \underline{\cal C}={\rm Im}(p)$ 
is an $A$-coring.   $(A,C,\psi)$ is a weak entwining structure if
and only if $A\otimes C$ is a pre-coring with structure maps defined as
in Proposition~\ref{prop.tak}. Proposition~\ref{weak} can be deduced
from this result. Since pre-corings appear naturally in the context of
weak entwining structures it would be interesting to 
study general properties of pre-corings along the lines of 
 the present paper or of \cite{Swe:pre}. 
These are the topics for 
further work.

\begin{center}
{\sc Acknowledgements}
\end{center}
I would like to thank  Stefaan Caenepeel for discussions and 
EPSRC for an Advanced Research Fellowship.

\end{document}